*Automated Reidemeister Moves: A Numerical Approach to the Unknotting Problem*

Dana Foley
The College of Wooster
1189 Beall Ave.
Wooster, OH, 44691
330-601-3368

**ABSTRACT**
In mathematics, a knot is a single strand of string crossed over itself any number of times, and connected at the ends. The Reidemeister Moves have been proven to be the three core moves necessary to fully untangle a knot. Some knots can be untangled to a loop (the unknot), while others are fundamentally knotted. We define four generalized moves based on the Reidemeister Moves. These moves have the capability of untangling knots which must be made more complicated before they can be simplified. With these moves, we construct a computer program which reads the two-dimensional projection of a knot in its Gauss Code notation and untangles it to the fewest possible number of crossings. Due to the properties of the Gauss Code notation, the program runs efficiently with minimal computation time, compared to currently existing untangling programs. We have tested it on all possible Gauss Codes of up to 50 crossings, and the program successfully determines a simplest possible projection of each knot. The results are consistent with our predictions about how the moves work in untangling a knot.

**1. INTRODUCTION**
One of the core problems in mathematical knot theory is determining whether or not a given knot is the unknot – a closed loop with no crossings. Ideally, we would like to have a computer program that can take a knot and untangle it to the fewest possible number of crossings. This problem can be, and has been, approached in many ways. Several existing programs such as the widely-used *KnotPlot* [4] consider a knot by its form in three-dimensional space and use the physical properties and interactions of the knot to untangle it. Wolfgang Haken proposed another method. Haken wrote an algorithm based on the so-called "bounding surface" of the knot. However, it is so complex that, in a knot whose projection has more than just 5 crossings, computer programs attempting to implement it are far too resource-intensive to be practical [2].

When untangling a mathematical knot, the typical method is to consider the knot's two-dimensional projection, and to use a series of what are called *Reidemeister moves*. This is easy and practical to do manually for knot projections with a low number of crossings, because we can visualize the two-dimensional knot projection in three-dimensional space and predict what moves we will need to perform in order to untangle the knot. Computer programs implementing the Reidemeister moves could easily handle simple knots when no moves need to occur which increase the number of crossings on the way to simplifying the knot. However, there are many cases in which such moves must occur. Knots that require these moves are called **Hard Unknots** [3]. This is where a computer program attempting to monotonically reduce the number of crossings is likely to fail, and where some others simply become too resource-intensive.

We introduce generalized versions of the Reidemeister Moves to make computer processing of a knot more practical. In this paper, we outline a program in Python that uses these modified Reidemeister Moves to efficiently and quickly untangle a knot. To implement this, we need to determine a method of representing knots in a more mathematical form.

## 2. GAUSS CODE – REPRESENTING A KNOT NUMERICALLY

Before our program can untangle a knot, we must decide how the computer will "see" knots. When untangling a knot by hand, we can visualize the process and notice patterns that allow us to choose the ideal moves to untangle the knot. We must determine how to numerically or algebraically represent knots so that the computer will recognize these patterns. We find that one of the simplest and most useful methods is the knot representation *Gauss Code*.

**Gauss Code** is a numerical way to represent a knot projection defined by the following procedure. To construct the Gauss Code of a knot, begin by picking a point on a two-dimensional knot projection, and a direction; both are arbitrary. Trace through the knot in the chosen direction. When a crossing is encountered that has not already been named, label it with a sequential integer beginning with 1. When encountering the overcrossing component of that crossing, denote it with the positive form of the corresponding integer. When encountering the undercrossing, instead write the negative number. Continue to trace through the knot, writing every positive or negative number in a list. When the first number in the list is reached again, stop tracing (and do not include the first number a second time). The list of numbers should contain a positive and negative number for every crossing, making the length of the list twice the number of crossings in the knot.

For example, Figure 1 demonstrates how we would trace through the trefoil knot.

The Gauss Code for the trefoil knot would be: 1,-2,3,-1,2,-3.

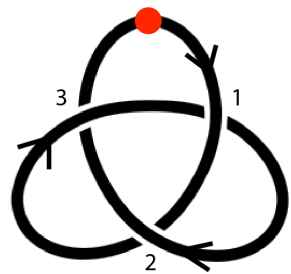

**Figure 1. Tracing through and labeling crossings of the trefoil knot, with a starting point at the top of the projection. This yields the Gauss Code: 1,-2,3,-1,2,-3**

## 3. STANDARD REIDEMEISTER MOVES IN GAUSS CODE

Certain patterns in the Gauss Code indicate where and how the Reidemeister Moves may be applied. We will examine each individually.

### 3.1 Move 1

The first Reidemeister Move is shown in Figure 2. If we were to denote this segment of the knot in Gauss Code, tracing from the bottom upward would yield $-a,a$, and from the top downward would be $a,-a$.

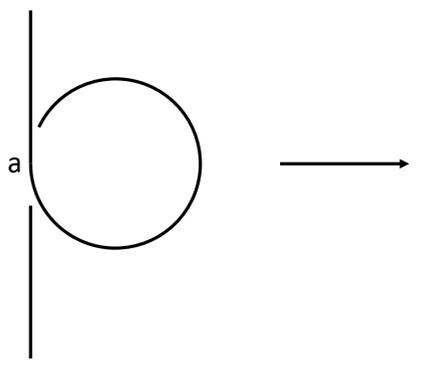

**Figure 2. Performing the first Reidemeister Move to remove a crossing**

This pattern, appearing anywhere in the Gauss Code, indicates that the first

Reidemeister Move may be performed. The result is the removal of the crossing *a*, thus deleting *a* and –*a* from the Gauss Code list.

To identify and perform the first Reidemeister Move, the program searches for two adjacent integers in the name which are negatives of each other. If this scenario is found, the numbers get removed.

### 3.2 Move 2

The second Reidemeister Move is shown in Figure 3. To denote this in Gauss Code, we must consider what would happen when tracing through each of the two strands - the vertical segment, and the curved one. The curved segment could be *a,b* or *b,a*, and the vertical segment could be –*a,*–*b*, or –*b,*–*a*. Any combination of these situations is valid. One can see readily that any appearance of this pattern in the Gauss Code indicates the opportunity to perform a second Reidemeister Move.

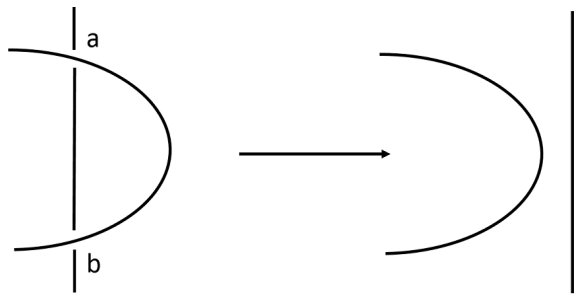

**Figure 3. Performing the second Reidemeister Move to remove two crossings**

To perform the second Reidemeister Move, the program looks for two adjacent crossings with the same sign (positive or negative). It then locates the negatives of these integers, and determines if those numbers are also adjacent in the list. If these conditions are true, it removes the four numbers involved from the list.

### 3.3 Move 3

The Third Reidemeister Move is shown in Figure 4.

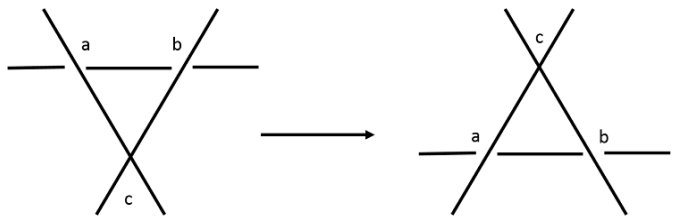

**Figure 4. Performing the third Reidemeister Move to relocate two crossings**

Note that the horizontal underpass passing through crossings *a* and *b* could instead be an overpass, and the form of crossing *c* is irrelevant. The third Reidemeister Move is handled as a particular type of case of our generalized moves, so we will not consider its direct implementation.

## 4. TANGLES

Tangles are the core of how the program will perform the generalized Reidemeister Moves, so we must define a tangle before we can explore those moves.

**Definition:** A **Tangle** is a region of a knot enclosed by a circle such that the knot crosses the circle exactly four times [5]. We define the **size** of a tangle as the number of crossings within the circle.

Figure 5 shows three examples of tangles, though infinitely many exist. We consider a tangle to be any interaction of the four entry strands within the circle.

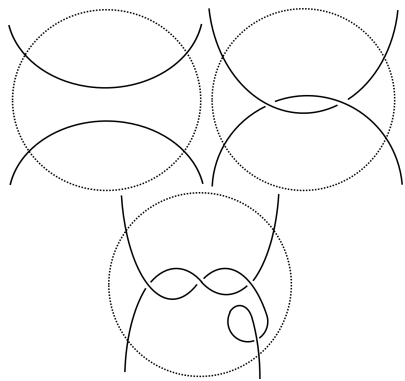

**Figure 5. Three simple tangles, with sizes 0, 2, and 4, respectively**

We have equipped the program with a function capable of locating every valid tangle in the knot, for every possible tangle size, as defined above. To do this, we scan every crossing in the knot. For each crossing, we trace forward over the specified size of the tangle in consideration. This sub-segment is the first of the two strands which compose the tangle. To find the second strand, we locate the first and last crossings in the first strand that are not completed within the same strand. That is to say, if a crossing $x$ is encountered in the first strand list, the crossing -$x$ is not in that list. Given these beginning and end points, we find the negative of each crossing number. We trace in the appropriate directions for a certain distance and list each crossing along the way. If the crossings in this second strand are all completed in either itself or the first strand, and the same is true for the first segment with respect to the second, the program declares this a valid tangle.

Within the cases in Figure 5, the first image is not recognized by the program as a tangle because it has no crossings. While technically it is a tangle by our definition, it will not be useful to consider when applying our generalized moves.

The second is a common, simple tangle. The crossings in each strand are all completed by the other strand, and not by itself. It is important to note that the number of crossings between the two strands is even, which results in the ends of each strand being directly adjacent.

The third tangle is somewhat different. The total number of crossings in the tangle is 4, an even number, but the number of crossings completed by the other strand is 3. Since this is an odd number, the ends of each strand are not adjacent. In this case, one strand has ends in the bottom left and top right of the circle, and the other strand has ends in the top left and bottom right.

## 5. GENERALIZED MOVE 1

Now, we propose two generalized translation moves that can replace the basic Reidemeister Moves in the process of untangling hard unknots, preventing the knot from ever becoming more complicated than it previously was.

The first Reidemeister Move remains as it is. However, we will now begin to consider it to be under a general "Move 1" category. The first Reidemeister Move reduces the number of crossings in a knot, so we will refer to it as **Reduction Move 1**.

We not propose a generalized move which alters the projection of a knot, but not its number of crossings. As such, we will call this move **Translation Move 1**.

Translation Move 1, as shown in Figure 6, takes a tangle which is adjacent to a crossing on one side, and flips the tangle over in the appropriate direction based on the crossing. This removes the original crossing and creates one on the opposite side of the tangle, essentially relocating the crossing across the tangle. The "flipping" nature of this move is why we define it to be in the "Move 1" category.

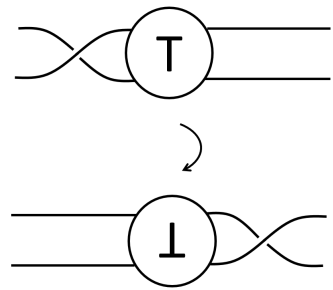

**Figure 6. Performing Translation Move 1 over the tangle**

To perform this computationally, the program looks for the crossing immediately beyond each end of a tangle. If any two of the crossing integers beyond the ends of the tangle are negations of each other, then we can perform a valid Translation Move 1. The program first negates the value of every number composing the tangle, flipping it over. Next, the two numbers composing the crossing adjacent to the tangle are removed from their current positions and reinserted in the appropriate locations on the opposite side of the tangle, thus completing the move as it would be seen in Figure 6.

## 6. GENERALIZED MOVE 2

Similarly, the second Reidemeister Move remains the same, and we refer to it as **Reduction Move 2**. We now propose the addition of **Translation Move 2**.

**Translation Move 2** is shown in Figure 7. To perform this move, we must consider a tangle and two of its adjacent ends. If the crossings immediately beyond these ends are both underpasses or overpasses, we can relocate the strand intersecting both tangle segments to the opposite side of the tangle. This act of relocating a strand with two adjacent over/underpasses, similar to Reduction Move 2, is why we choose to place Translation Move 2 in the same category.

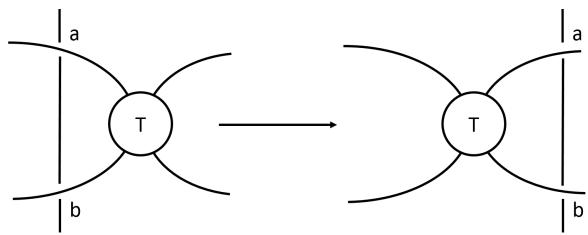

**Figure 7. Performing Translation Move 2 over the tangle *T***

To perform this move computationally, we again examine what crossings occur immediately beyond two adjacent ends of the tangle. If the two crossings are both overpasses or underpasses (both positive or negative), then we check the negative of each number, and determine if those crossings are adjacent. This process should seem similar to how we check the conditions for Reduction Move 2. If a Translation Move 2 is available, performing it consists of relocating the two crossings composing the overpass/underpass to the other side of the knot. In the case of Figure 7, we would relocate the crossings *a* and *b* to the other side of the tangle. However, the numbers *-a* and *-b* remain as they are, because moving this underpass strand does not affect how *-a* and *-b* are read in the Gauss Code of the knot.

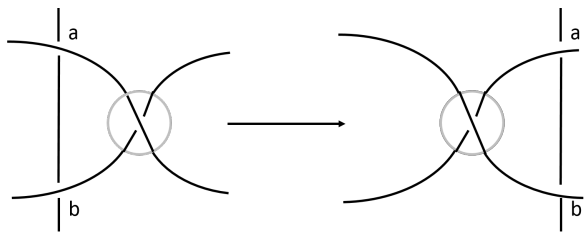

**Figure 8. A particular case of Translation Move 2 where *T* is a single crossing, demonstrating the presence of the original third Reidemeister Move**

Referring back to the original third Reidemeister Move, it turns out to be included in the general Translation Move 2. Allowing the tangle *T* from Figure 7 to be a single crossing results in the scenario seen in Figure 8. This is the appropriate form of the third Reidemeister Move, and our method for performing Translation Move 2 allows this to work as intended.

## 7. PROGRAM STRUCTURE

Now that our program is able to recognize and perform the four moves we have defined, we need to set up a structure in Python for how and when it will perform each move.

Our program accepts the input for a knot in the form of the Gauss Code list of numbers. First, it runs a check process to make sure the list "makes sense"; all crossing numbers must have exactly one corresponding negative number that completes the crossing. If the knot is not valid, the program terminates.

If the knot is valid, we move on to a While loop which repeats while the knot is still tangled (not the unknot) or while the program has not tried every untangling operation it is capable of.

In this While loop, we are repeating a "step" function, which checks for and performs any moves necessary to reduce the knot. If the step function fails, the knot has been untangled as much as the program is capable of, and theoretically to the fewest possible number of crossings. When the step function fails, the program terminates. This is an advantage over Haken's algorithm [2], which cannot recognize knots other than the unknot.

Within the step function, several checks run in a specific order. First, it checks if a Reduction Move 2 can be performed. If so, it performs the move and ends the function, returning to the While loop. If Reduction Move 2 fails, it does the same check for a Reduction Move 1. We try Reduction Move 2 first because it removes more crossings per step, slightly reducing run time.

If no reduction moves are available, it runs through a complex set of steps which examines possible translation moves. First, it tries every possible translation move (1 and 2) independently. If performing any of those moves results in a reduction move becoming available, that translation move is performed, and the function ends, returning to the While loop so that the reduction move may be performed next.

There are rare cases where, to perform a reduction move, multiple instances of translation moves must be performed. These situations are uncommon, so the program usually does not reach this particular computation at any point when untangling the unknot. However, if needed, the program will examine every possible projection of the knot with the current number of crossings, reachable through instances of Translation Moves 1 and 2. If none of these allow a reduction move to occur, the program declares the knot fully untangled, and the while loop terminates. The program then returns the final knot in Gauss Code.

## 8. RESULTS AND FUTURE WORK

We have tested this program on every possible Gauss Code for up to 50 crossings. In each case, the program successfully reduced every knot to the least possible number of crossings. This suggests that our program could potentially be a solution to the unknotting problem, especially considering the great number of hard unknots present in the list of knots the program has tested.

In most cases, the step function of our program only requires one translation move before a reduction move can occur. In the cases where more than one translation move is required, it has (so far) always been a sequence of Translation Move 2s, where a tangle was brought across multiple over/underpasses in the same direction repeatedly. For example, a tangle would need to be moved across an overpass on the left, then across an underpass on the left, etc. We speculate that other combinations of translation moves may not be necessary, but the full check remains a part of the program, to be more certain that our results are accurate.

The extra check only affects computation time slightly. The average time to fully untangle several unknot diagrams of 50 crossings was approximately 0.3 seconds for each knot. The program is equipped to read a large (>1000) list of knots and untangle them all, with computation time of approximately 30 seconds, depending on the complexity of each knot. This is quite a significant advantage over existing implementations of Haken's algorithm, which take around 3 days to identify the unknot from a knot with a single crossing [2].

In the future, we hope to expand upon our predictions about combinations of translation moves to further optimize the program. Even more ideally, we hope to determine a way to prove that Translation Moves 1 and 2 cover all possible scenarios where a knot must be made more complicated before it can be simplified. This would, in theory, be a much stronger proposal that our computer program solves the unknotting problem.